\documentclass[12pt,reqno]{amsart}

\usepackage{graphicx}
\graphicspath{ {./images/} }
\usepackage{bm}

\usepackage{color}

\usepackage{amssymb}

\usepackage{amsfonts}

\usepackage{amsmath}

\usepackage{amsthm}

\usepackage[all]{xy}

\usepackage[mathlines]{lineno}

\usepackage{hyperref}

\newtheorem{theorem}{Theorem}[section]

\newtheorem{corollary}[theorem]{Corollary}

\newtheorem{lemma}[theorem]{Lemma}

\newtheorem{proposition}[theorem]{Proposition}

\newtheorem{conjecture}[theorem]{Conjecture}

\newtheorem{Definition}[theorem]{Definition}

\newtheorem{Example}[theorem]{Example}

\newtheorem{Remark}[theorem]{Remark}

\newenvironment{remark}{\begin{Remark}\begin{em}}{\end{em}\end{Remark}}

\def\diag{{\mbox{diag}\,}}

\DeclareMathOperator{\tr}{tr}

\begin{document}

\title[Weak log-majorization among three means]{Weak log-majorization between the geometric and Wasserstein means}

\author{Luyining Gan}
\address{Luyining Gan, School of Science, Beijing University of Posts and Telecommunications, Beijing 100876, China \& Department of Mathematics and Statistics, University of Nevada at Reno, Reno, NV 89557, USA}
\email{luyining.gan@gmail.com; lgan@unr.edu}

\author{Sejong Kim}
\address{Sejong Kim, Department of Mathematics, Chungbuk National University, Cheongju 28644, Korea}
\email{skim@chungbuk.ac.kr}
\maketitle

\begin{abstract}
There exist lots of distinct geometric means on the cone of positive definite Hermitian matrices such as the metric geometric mean, spectral geometric mean, log-Euclidean mean and Wasserstein mean.
In this paper, we prove the log-majorization relation on the singular values of the product of given two positive definite matrices and their (metric and spectral) geometric means. We also establish the weak log-majorization between the spectra of two-variable Wasserstein mean and spectral geometric mean. In particular, we verify with certain condition on variables that
two-variable Wasserstein mean converges decreasingly to the log-Euclidean mean with respect to the weak log-majorization.
\end{abstract}



\section{Introduction}

Let $\mathbb{C}_{m \times m}$ be the space of all $m \times m$ complex matrices,
and $\mathbb{H}_m$ be the real vector space of all $m \times m$ Hermitian matrices.
We denote as $\mathbb{P}_m \subset \mathbb{H}_m$ the open convex cone of all $m \times m$ positive definite matrices.
Given $A \in \mathbb{H}_{m}$, we use $A \geq (>) 0$ to indicate that $A$ is positive semidefinite (positive definite, respectively).
For $A, B \in \mathbb{H}_m$, the Loewner order $A \geq B$ means $A-B \geq 0$, that is, $A-B$ is positive semidefinite.

For $X \in \mathbb{C}_{m \times m}$, the singular values of $X$ are the eigenvalues of $|X| := (X^*X)^{1/2}$.
Denote by $s(X)$ the $m$-tuple of all singular values of $X \in \mathbb{C}_{m \times m}$ with non-increasing order: $s_1(X) \geq s_2(X) \geq \cdots \geq s_m(X) \geq 0$. We also denote by $\lambda(X)$ the $m$-tuple of all real eigenvalues of $X \in \mathbb{H}_{m}$ with $\lambda_1(X) \geq \lambda_2(X) \geq \cdots \geq \lambda_m(X)$.

Let us recall the definition of (weak) log-majorization.
Let $x, y$ be two $m$-tuples of positive real numbers. Denote by $x^{\downarrow}, y^{\downarrow}$ the non-increasing order of elements of $x, y$ respectively. We write $x \prec_{w \log} y$ if $x$ is \emph{weakly log-majorized} by $y$, that is,
\begin{equation}\label{eqn:wlog}
\prod_{i=1}^k x^{\downarrow}_i \leq \prod_{i=1}^k y^{\downarrow}_i, \quad k = 1, 2, \dots, m.
\end{equation}
We say that $x$ is \emph{log-majorized} by $y$, denoted by $x \prec_{\log} y$, if \eqref{eqn:wlog} is true for $k = 1, 2, \dots, m-1$ and equality holds for $k=m$.
For simplicity, we write $A \prec_{\log} B$ if $\lambda(A) \prec_{ \log} \lambda(B)$, and $A \prec_{w \log} B$ if $\lambda(A) \prec_{w \log} \lambda(B)$ for $A, B \in \mathbb{P}_m$.

Let $A, B \in \mathbb{P}_{m}$ and $t \in \mathbb{R}$. The \emph{metric geometric mean} of $A, B$ is a differentiable curve on $\mathbb{P}_m$ defined by
\[ A \#_t B = A^{1/2} (A^{-1/2} B A^{-1/2})^t A^{1/2}.
\]
This notion was first introduced by Pusz and Woronowicz \cite{PW75} for $t = 1/2$, simply denoted as $A \# B = A \#_{1/2} B$.
The weighted version was introduced later by Kubo and Ando \cite{KA79}.
As another notion of geometric mean on $\mathbb{P}_m$, the \emph{spectral geometric mean} of $A, B$ is a differentiable curve defined by
\[ A \natural_t B = (A^{-1} \# B)^t A (A^{-1} \# B)^t,
\]
which was first proposed by Fiedler and Pt\'ak \cite{FP97} for the version of $t = 1/2$. We simply denote as $A \natural B = A \natural_{1/2} B$. The weighted version was introduced later by Lee and Lim \cite{HL07}, and its several properties have been recently established \cite{GK23, Kim21} on the setting of positive invertible operators.

Note that the metric and spectral geometric means can be considered as the non-commutative versions for geometric mean of positive scalars. Zou \cite{ZZ17} provided their relationship in terms of the log-majorization on the singular values of the metric geometric mean and the multiplication of two matrices, that is, for any $A, B \in \mathbb{P}_{m}$,
\[
s(A^{1/2} (A \# B) B^{1/2}) \prec_{\log} s(A B).
\]
It is straightforward to consider such relation to the weighted version and to extend such relation to the spectral geometric mean.

The \emph{Wasserstein distance} of $A, B \in \mathbb{P}_m$ is the Riemannian metric given by
\begin{displaymath}
d_{W}(A, B) = \left[ \tr \left( \frac{A + B}{2} \right) - \tr (A^{1/2} B A^{1/2})^{1/2} \right]^{1/2}.
\end{displaymath}
This coincides with the Bures distance of density matrices in quantum information theory, and can be considered as a matrix version of the Hellinger distance for probability vectors.
The Wasserstein mean of $A_{1}, \dots, A_{n} \in \mathbb{P}_{m}$ is the least squares mean for the Wasserstein distance, defined by
\begin{displaymath}
\Omega(\omega; A_{1}, \dots, A_{n}) = \underset{X \in \mathbb{P}_{m}}{\arg \min} \sum_{j=1}^{n} w_{j} d_{W}^{2}(X, A_{j}),
\end{displaymath}
where $\omega = (w_{1}, \dots, w_{n})$ is a positive probability vector.
In particular, when $n=2$, replacing $A_1, A_2$ with $A, B$ and $(w_1, w_2)$ with $(1-t, t)$, two-variable \emph{Wasserstein mean} of $A$ and $B$ has the following explicit formula:
\begin{displaymath}
A \diamond_{t} B := \Omega((1-t, t); A, B) = (1-t)^{2} A + t^{2} B + t(1-t) [A (A^{-1} \# B) + (A^{-1} \# B) A].
\end{displaymath}
The numerical computation and applications of the Wasserstein mean from both theoretical and computational aspects have been widely studied: see~\cite{ADCM16, BJL19, BJL19a, CD14, HK19, ZCY22} and references therein.

Recently, the following (weak) log-majorization relation among matrix means has been shown \cite{BJL19}:
\begin{displaymath}
A \#_{t} B \prec_{\log} \exp ((1-t) \log A + t \log B) \prec_{w \log} A \diamond_{t} B.
\end{displaymath}
Moreover, the monotonicity of metric and spectral geometric means with respect to the log-majorization has been proved \cite{AH94, GLT21}:
\begin{displaymath}
A \#_{t} B \nearrow_{\prec_{\log}} \exp ((1-t) \log A + t \log B) \swarrow_{\prec_{\log}} A \natural_{t} B.
\end{displaymath}
So finding a weak log-majorization relation between the spectral geometric and Wasserstein mean, and proving the monotonicity of Wasserstein mean with respect to the weak log-majorization are interesting problems.

In this paper, we organize the sections as follows.
In Section~\ref{Sec:properties}, we recall known results for the metric geometric mean, spectral geometric mean and Wasserstein mean of positive definite matrices with log-majorization and fundamental properties.
We then prove the log-majorization relation on the singular values of geometric means and the product of matrices in Section~\ref{Sec:geomeans}. The main goal is to establish the weak log-majorization relation between the spectral geometric mean and Wasserstein mean in Section \ref{Sec:log_maj} and the monotonicity of Wasserstein mean with respect to the weak log-majorization under certain condition in Section~\ref{Sec:Was}.


\section{Preliminaries on log-majorization and matrix means} \label{Sec:properties}

There has been a focus on the study of the (weak) log-majorization relations between different means on $\mathbb{P}_m$.
We know from \cite{BJL19} the following (weak) log-majorization relation among the Cartan (Riemannian) mean $\Lambda$, log-Euclidean mean $L$ and Wasserstein mean $\Omega$ of $A_1, \dots, A_n \in \mathbb{P}_m$:
\begin{equation} \label{eqn:Cartan-Log-Wass}
\Lambda(\omega; A_1, \dots, A_n) \prec_{\log} L(\omega; A_1, \dots, A_n) \prec_{w \log} \Omega(\omega; A_1, \dots, A_n),
\end{equation}
where
\begin{displaymath}
\Lambda(\omega; A_1, \dots, A_n) := \underset{X \in \mathbb{P}_{m}}{\arg \min} \sum_{j=1}^{n} w_{j} d_{R}^{2}(X, A_{j})
\end{displaymath}
is the Cartan mean for the Riemannian trace metric $d_{R}(A, B) = \Vert \log A^{-1/2} B A^{-1/2} \Vert_{2}$ and
\begin{displaymath}
L(\omega; A_1, \dots, A_n) := \exp \left( \sum_{j=1}^{n} w_{j} \log A_{j} \right)
\end{displaymath}
is the log-Euclidean mean. In particular, replacing $A_1, A_2$ with $A, B$ and $(w_1, w_2)$ with $(1-t, t)$ in \eqref{eqn:Cartan-Log-Wass} for $n = 2$, we obtain
\begin{equation*} \label{eqn:g-log-Wass}
A \#_{t} B \prec_{\log} \exp((1-t) \log A + t \log B) \prec_{w \log} A \diamond_t B.
\end{equation*}
Furthermore, the log-majorization between the log-Euclidean mean and spectral geometric mean has been shown in \cite{GT22}:
\begin{displaymath}
\exp((1-t) \log A + t \log B) \prec_{\log} A \natural_t B.
\end{displaymath}
So it is a natural question whether there exists a weak log-majorization relation between the spectral geometric mean and the Wasserstein mean.

We see some properties of the metric geometric mean \cite{Bh07}, spectral geometric mean \cite{Kim21, HL07}, and Wasserstein mean \cite{HK19, HK22}, which are useful to prove our main results.

\begin{lemma} \label{L:m-geo}
Let $A, B, C, D \in \mathbb{P}_{m}$ and let $s, t, u \in [0,1]$. Then the following are satisfied.
\begin{enumerate}
\item $A \#_t B = B \#_{1-t} A$.
\item $(A \#_t B)^{-1} = A^{-1} \#_t B^{-1}$.
\item $A \#_t B \leq C \#_t D$ whenever $A \leq C$ and $B \leq D$.
\item $M (A \#_t B) M^* = (M A M^*) \#_t (M B M^*)$ for any non-singular matrix $M$.
\item $(aA) \#_t (bB) = a^{1-t} b^{t} (A \#_t B)$ for any $a, b >0$.
\item $(A \#_s B) \#_t (A \#_u B) = A \#_{(1-t)s+tu} B$.
\item $\det (A \#_{t} B) = (\det A)^{1-t} (\det B)^{t}$.
\item $((1-t)A^{-1}+tB^{-1})^{-1}\leq A\#_t B \leq (1-t)A +tB$
\end{enumerate}
\end{lemma}

\begin{lemma} \label{L:s-geo}
Let $A, B \in \mathbb{P}_{m}$ and let $s, t, u \in [0,1]$. Then the following are satisfied.
\begin{enumerate}
\item $A \natural_t B = B \natural_{1-t} A$.
\item $(A \natural_t B)^{-1} = A^{-1} \natural_t B^{-1}$.
\item $U (A \natural_t B) U^* = (U A U^*) \natural_t (U B U^*)$ for any unitary matrix $U$.
\item $(aA) \natural_t (bB) = a^{1-t} b^{t} (A \natural_t B)$ for any $a, b >0$.
\item $(A \natural_s B) \natural_t (A \natural_u B) = A \natural_{(1-t)s+tu} B$.
\item $\det (A \natural_{t} B) = (\det A)^{1-t} (\det B)^{t}$.
\end{enumerate}
\end{lemma}

\begin{lemma} \label{L:Wass}
Let $A, B \in \mathbb{P}_m$ and let $s, t, u \in [0,1]$. Then the following are satisfied.
\begin{enumerate}
\item $A \diamond_t B = B \diamond_{1-t} A$.
\item $(A \diamond_t B)^{-1} = A^{-1} \diamond_t B^{-1}$ if and only if $A = B$.
\item $U (A \diamond_t B) U^* = (U A U^*) \diamond_t (U B U^*)$ for any unitary matrix $U$.
\item $(aA)\diamond_t (aB) = a (A \diamond_t B)$ for any $a > 0$.
\item $(A \diamond_s B) \diamond_t (A \diamond_u B) = A \diamond_{(1-t)s+tu} B$.
\item $\det (A \diamond_{t} B) \geq (\det A)^{1-t} (\det B)^{t}$.
\end{enumerate}
\end{lemma}


\section{Log-majorization of geometric means} \label{Sec:geomeans}

Note that the metric geometric mean $A \#_{t} B$ and spectral geometric mean $A \natural_{t} B$ are non-commutative versions for geometric mean of positive scalars. In other words, for commuting $A, B \in \mathbb{P}_{m}$ the metric geometric mean and spectral geometric mean become $A^{1-t} B^{t}$. So it is an interesting problem to find the relationship between the non-commutative and commutative versions of geometric mean.

Zou \cite{ZZ17} proved that for any $A, B \geq 0$,
\begin{equation} \label{eqn:Zou}
s(A^{1/2} (A \# B) B^{1/2}) \prec_{\log} s(A B).
\end{equation}
Lemos and Soares \cite{LS18} provided another proof of \eqref{eqn:Zou}, and more generally asked whether there exists the following log-majorization relation for $A, B \geq 0$
\begin{equation} \label{ineq:LS}
s (A^{t} (A \#_{t} B) B^{1-t}) \prec_{\log} s(A B),  \quad t \in [0,1].
\end{equation}
This is still an open problem, but Ghabries et al. \cite{GA22} established the following log-majorization results related to \eqref{ineq:LS}.
\begin{theorem}
Let $A, B \geq 0$.
\begin{enumerate}
\item $s(A^t (A \#_t B) B^{1-t}) \prec_{\log} s(A^{\frac 3 2} B A^{-\frac 1 2})\quad \text{for} \quad \frac 1 2 \leq t \leq 1$,
\item $s(A^t (A \#_t B) B^{1-t}) \prec_{\log} s(B^{\frac 3 2} A B^{-\frac 1 2})\quad \text{for} \quad 0 \leq t \leq \frac 1 2$.
\end{enumerate}
\end{theorem}

We prove alternative versions of \eqref{ineq:LS} for the metric geometric and spectral geometric means. To prove this, we introduce the following well-known inequality.

\begin{theorem}[Loewner-Heinz inequality]\label{thm:LH-ineq}
Let $A \geq B\geq 0$ and $p \in[0, 1]$. Then
$A^p \geq B^p$.
\end{theorem}

Antisymmetric tensor power is a standard technique in the theory of log-majorization. For the positive semidefinite matrices, there are some interesting properties. Note that for $A\geq 0$ and $1\leq k \leq m$,
\[ \prod_{i=1}^k \lambda_i(A) = \lambda_1(\Lambda^k A),
\]
where $\Lambda^k A$ is the $k$th antisymmetric tensor power (or the $k$th compound matrix) of $A$.
By the definition of log-majorization, for $A, B\geq 0$, $A\prec_{\log} B$ if and only if $\lambda_1(\Lambda^k A) \leq \lambda_1(\Lambda^k B)$, $k = 1, . . . , m-1$, and $\det A=\det B$.
Moreover, the map $\mathbb{P}_{m} \ni A \mapsto \Lambda^k A$ is multiplicative, that is,
\[ \Lambda^k(AB) = \Lambda^k(A) \Lambda^k(B) \quad \text{and} \quad (\Lambda^k(A))^r = \Lambda^k(A^r), \ r \in (-\infty, \infty).
\]
So it is clear that $\Lambda^k(A\#_t B) = (\Lambda^k(A))\#_t (\Lambda^k(B))$ and $\Lambda^k(A\natural_t B) = (\Lambda^k(A))\natural_t (\Lambda^k(B))$ for $A, B\in \mathbb{P}_m$.
Now we are ready to prove alternative versions of \eqref{ineq:LS} for the metric geometric and spectral geometric means.

\begin{theorem}
Let $A, B \in \mathbb{P}_{m}$. For any $t \in [0,1]$,
\begin{equation} \label{eqn:g-majorization}
s(A^{t - \frac{1}{2}} (A \#_{t} B) B^{\frac{1}{2} - t}) \prec_{\log} s(A^{1/2} B^{1/2}).
\end{equation}
\end{theorem}

\begin{proof}
Note that \eqref{eqn:g-majorization} is equivalent to
\begin{equation} \label{E:lambda-g-major}
\lambda(A^{t - \frac{1}{2}} (A \#_{t} B) B^{1-2t} (A \#_{t} B) A^{t - \frac{1}{2}})^{1/2} \prec_{\log} \lambda(A^{1/2} B A^{1/2})^{1/2}.
\end{equation}
Since $(A^{t - \frac{1}{2}} (A \#_{t} B) B^{1-2t} (A \#_{t} B) A^{t - \frac{1}{2}})^{1/2}$ and $(A^{1/2} B A^{1/2})^{1/2}$ are both homogeneous from Lemma \ref{L:m-geo} (5) and the $k$th antisymmetric tensor power preserves the matrices from two sides, it is enough to show that
\begin{center}
$A^{1/2} B A^{1/2} \leq I \qquad$ implies $\qquad A^{t - \frac{1}{2}} (A \#_{t} B) B^{1-2t} (A \#_{t} B) A^{t - \frac{1}{2}} \leq I$.
\end{center}

\begin{itemize}
\item[Step 1.] We first prove \eqref{E:lambda-g-major} for $t \in [0, 1/2]$. Assume that $A^{1/2} B A^{1/2} \leq I$. Then $B \leq A^{-1}$, and $B^{1-2t} \leq A^{2t-1}$ by Theorem~\ref{thm:LH-ineq} since $2t \in [0,1]$.
\begin{displaymath}
\begin{split}
A^{t - \frac{1}{2}} (A \#_{t} B) B^{1-2t} (A \#_{t} B) A^{t - \frac{1}{2}} & \leq A^{t - \frac{1}{2}} (A \#_{t} B) A^{2t-1} (A \#_{t} B) A^{t - \frac{1}{2}} \\
& = \left( A^{t - \frac{1}{2}} (A \#_{t} B) A^{t - \frac{1}{2}} \right)^{2}.
\end{split}
\end{displaymath}
Since $B \leq A^{-1}$, we obtain from Lemma \ref{L:m-geo} (3)
\begin{displaymath}
A^{t - \frac{1}{2}} (A \#_{t} B) A^{t - \frac{1}{2}} \leq A^{t - \frac{1}{2}} A^{1-2t} A^{t - \frac{1}{2}} \leq I.
\end{displaymath}
Therefore, $A^{t - \frac{1}{2}} (A \#_{t} B) B^{1 - 2t} (A \#_{t} B) A^{t - \frac{1}{2}} \leq \left( A^{t - \frac{1}{2}} (A \#_{t} B) A^{t - \frac{1}{2}} \right)^{2} \leq I$. Moreover, by Lemma \ref{L:m-geo} (7)
\begin{displaymath}
\det \left[ A^{t - \frac{1}{2}} (A \#_{t} B) B^{1-2t} (A \#_{t} B) A^{t - \frac{1}{2}} \right] = \det (A B) = \det (A^{1/2} B A^{1/2}),
\end{displaymath}
and thus, \eqref{E:lambda-g-major} holds for $t \in [0, 1/2]$.

\item[Step 2.] Let $t \in [1/2, 1]$. Since $s(X) = s(X^{*})$ for any matrix $X \in \mathbb{C}_{m \times m}$, we have from Lemma \ref{L:m-geo} (1) and Step 1
\begin{displaymath}
\begin{split}
s(A^{t - \frac{1}{2}} (A \#_{t} B) B^{\frac{1}{2} - t})  = s(B^{\frac{1}{2} - t} (B \#_{1-t} A) A^{t - \frac{1}{2}}) & = s(B^{(1-t)-\frac{1}{2}} (B \#_{1-t} A) A^{ \frac{1}{2}- (1-t)})\\
& \prec_{\log} s(B^{1/2} A^{1/2}) = s(A^{1/2} B^{1/2}),
\end{split}
\end{displaymath}
which completes the proof.\qedhere
\end{itemize}
\end{proof}

\begin{lemma} \label{L:natural_ineq}
Let $A, B \in \mathbb{P}_m$. If $A \leq I$ and $B \leq I$, then $A \natural B \leq I$.
\end{lemma}

\begin{proof}
Let $A \leq I$ and $B \leq I$. Since the square of $A \natural B$ is similar to $A B$ by \cite{FP97},
\[ \lambda_1(A \natural B) = \lambda_1^{1/2}(A B) \leq \lambda_1^{1/2}(A) \lambda_1^{1/2}(B) \leq 1,\]
where the first inequality follows from $\lambda_{1}(A B) \leq \lambda_{1}(A) \lambda_{1}(B)$ for $A, B \in \mathbb{P}_{m}$. This implies that $A \natural B \leq I$.
\end{proof}

\begin{theorem} \label{T:s-geomean-major}
Let $A, B\in \mathbb{P}_{m}$ with $A \geq I$. Let $1/2 \leq t \leq 1$ and $0 \leq u \leq 1/2$. Then
\begin{equation} \label{eqn:s-majorization}
s(A^{-u} (A \natural_{t} B) B^{u}) \prec_{w \log} s(A^{1/2} B^{1/2}).
\end{equation}
In addition, if $\det A = \det B$ or $t = 1/2, u = 0$ then
\begin{displaymath}
s(A^{-u} (A \natural_{t} B) B^{u}) \prec_{\log} s(A^{1/2} B^{1/2}).
\end{displaymath}
\end{theorem}

\begin{proof}
Note that \eqref{eqn:s-majorization} is equivalent to
\begin{equation*} \label{E:lambda-s-major}
\lambda(A^{-u} (A \natural_{t} B) B^{2u} (A \natural_{t} B) A^{-u})^{1/2} \prec_{\log} \lambda(A^{1/2} B A^{1/2})^{1/2}.
\end{equation*}
Since $(A^{-u} (A \natural_{t} B) B^{2u} (A \natural_{t} B) A^{-u})^{1/2}$ and $(A^{1/2} B A^{1/2})^{1/2}$ are both homogeneous by Lemma \ref{L:s-geo} (4), it is enough to show that
\begin{center}
$A^{1/2} B A^{1/2} \leq I \qquad$ implies $\qquad A^{-u} (A \natural_{t} B) B^{2u} (A \natural_{t} B) A^{-u} \leq I$.
\end{center}

Let $A, B \in \mathbb{P}_{m}$ with $A \geq I$. Assume that $A^{1/2} B A^{1/2} \leq I$. Then $B \leq A^{-1}$, and $B^{2u} \leq A^{-2u}$ by Theorem~\ref{thm:LH-ineq} because $2u \in [0,1]$. So
\begin{displaymath}
A^{-u} (A \natural_{t} B) B^{2u} (A \natural_{t} B) A^{-u} \leq A^{-u} (A \natural_{t} B) A^{-2u} (A \natural_{t} B) A^{-u} = (A^{-u} (A \natural_{t} B) A^{-u})^{2}.
\end{displaymath}

Set $T := \{ t \in [\frac 12,1]: A \natural_{t} B \leq I \}$. Since $B \leq A^{-1}$ is equivalent to $A \natural_{1/2} B \leq I$ from \cite[Theorem 5]{Kim21}, we have $1/2 \in T$, and $1 \in T$ because $A \natural_{1} B = B \leq A^{-1} \leq I$. Assume that $s, t \in T$. Then by Lemma \ref{L:s-geo} (5) and Lemma \ref{L:natural_ineq}
\begin{displaymath}
A \natural_{\frac{s+t}{2}} B = (A \natural_{s} B) \natural_{1/2} (A \natural_{t} B) \leq I,
\end{displaymath}
so $\frac{s+t}{2} \in T$. This yields that $T$ contains all dyadic rational numbers in $[1/2, 1]$, and by the density of dyadic rational numbers and the continuity of spectral geometric mean $T = [1/2, 1]$.

Since $A \natural_{t} B \leq I$ for $1/2 \leq t \leq 1$, we obtain
\begin{displaymath}
A^{-u} (A \natural_{t} B) A^{-u} \leq A^{-2u} \leq I.
\end{displaymath}
Therefore, $A^{-u} (A \natural_{t} B) B^{2u} (A \natural_{t} B) A^{-u} \leq (A^{-u} (A \natural_{t} B) A^{-u})^{2} \leq I$. Moreover, in order that
\begin{displaymath}
\det (A^{-u} (A \natural_{t} B) B^{2u} (A \natural_{t} B) A^{-u})^{1/2} = (\det A)^{1-t-u} (\det B)^{t+u} = (\det A \det B)^{1/2},
\end{displaymath}
we have that $\det (A^{-1} B)^{t + u - 1/2} = 1$. So $\det A = \det B$, otherwise $t + u = 1/2$. Since $1/2 \leq t \leq 1$ and $0 \leq u \leq 1/2$, we obtain $t = 1/2$ and $u = 0$ if $\det A \neq \det B$.
\end{proof}

\begin{corollary}
Let $A, B \in \mathbb{P}_{m}$ with $B \geq I$. Let $0 \leq t \leq 1/2$ and $-1/2 \leq u \leq 0$. Then
\begin{equation*} \label{eqn:s-majorization-B}
s(A^{-u} (A \natural_{t} B) B^{u}) \prec_{w \log} s(A^{1/2} B^{1/2}).
\end{equation*}
\end{corollary}

\begin{proof}
By Lemma \ref{L:s-geo} (1) and Theorem \ref{T:s-geomean-major} with $B \geq I$, $1/2 \leq 1-t \leq 1$ and $0 \leq -u \leq 1/2$ we obtain
\begin{displaymath}
\begin{split}
s(A^{-u} (A \natural_{t} B) B^{u}) & = s(B^{u} (B \natural_{1-t} A) A^{-u}) \\
& \prec_{w \log} s(B^{1/2} A^{1/2}) = s(A^{1/2} B^{1/2}).\qedhere
\end{split}
\end{displaymath}
\end{proof}


\section{Weak log-majorization between two means} \label{Sec:log_maj}

We study in this section the relationship between the spectral geometric mean and the Wasserstein mean.
There are several different expressions of two-variable Wasserstein mean of positive invertible operators: see \cite[Lemma 2.4]{HK22} . We use the equivalent expressions as follows: for $A, B \in \mathbb{P}_m$ and $t \in [0,1]$
\begin{eqnarray}
A \diamond_{t} B &=& A^{-1/2} \left[ (1-t) A + t (A^{1/2} B A^{1/2})^{1/2} \right]^2 A^{-1/2} \label{eqn:Wass} \\
&=& [I \nabla_{t} (A^{-1} \# B)] A [I \nabla_{t} (A^{-1} \# B)] \label{eqn:Wass-v}
\end{eqnarray}
where $A \nabla_{t} B = (1-t) A + t B$ is the weighted arithmetic mean of $A$ and $B$. Note that the arithmetic-Wasserstein mean inequality has been shown in \cite{BJL19}:
\begin{equation} \label{E:A-Wass}
A \diamond_{t} B \leq A \nabla_{t} B.
\end{equation}

\begin{proposition} \label{P:sp-Wass}
Let $A, B \in \mathbb{P}_m$ and $t \in [0,1]$. If $A \diamond_{t} B \leq I$, then $A \natural_t B \leq I$.
\end{proposition}

\begin{proof}
For $t = 0$ and $t = 1$ it is obvious. We first consider the case $t \in (0, \frac{1}{2}]$.
Let $C := A^{-1} \# B$.
Suppose that $A \diamond_t B \leq I$. Note that $A^{1/2} C A^{1/2} = (A^{1/2} B A^{1/2})^{1/2}$.
Then by~\eqref{eqn:Wass}, Theorem~\ref{thm:LH-ineq} and taking the congruence transformation by $A^{-1/2}$, we obtain the following:
\begin{eqnarray*}
\left[ (1-t) A + t (A^{1/2} B A^{1/2})^{1/2} \right]^2 & \leq & A \\
(1-t) A + t (A^{1/2} B A^{1/2})^{1/2} & \leq & A^{1/2}  \\
(1-t) A + t A^{1/2} C A^{1/2} & \leq & A^{1/2} \\
(1-t) I + t C & \leq & A^{-1/2} \\
C & \leq & \frac{1}{t} A^{-1/2} + \left( 1 - \frac{1}{t} \right) I.
\end{eqnarray*}
Then we consider the largest eigenvalue of $A \natural_t B$. Since $2t \in (0,1]$, it can be computed as
\begin{equation} \label{eqn:eigspe}
\lambda_1(A \natural_t B) = \lambda_1(C^t A C^t) = \lambda_1(A^{1/2} C^{2t} A^{1/2}) \leq  \lambda_1 \left( \left(\frac{1}{t} A^{-1/2} + \left( 1 - \frac{1}{t} \right) I \right)^{2t} A \right).
\end{equation}
Since $A \in \mathbb{P}_m$, there exist a unitary matrix $U$ and a diagonal matrix $D = \diag (\lambda_1, \dots, \lambda_m)$ such that $A = U D U^*$. Then
\begin{eqnarray*}
& {} & \left( \frac{1}{t} A^{-1/2} + \left( 1 - \frac{1}{t} \right) I \right)^{2t} A \\
& = & U \left( \frac{1}{t} D^{-1/2} + \left( 1 - \frac{1}{t} \right) I \right)^{2t} D U^* \\
& = & U \begin{bmatrix} \left( \frac{1}{t} \lambda_1^{-1/2} + \left( 1 - \frac{1}{t} \right) \right)^{2t} \lambda_1 & & \\
& \ddots & \\
& & \left( \frac{1}{t} \lambda_m^{-1/2} + \left( 1 - \frac{1}{t} \right) \right)^{2t} \lambda_m \\
\end{bmatrix} U^*.
\end{eqnarray*}
We claim that $\displaystyle \left( \frac{1}{t} \lambda^{-1/2} + \left( 1 - \frac{1}{t} \right) \right)^{2t} \lambda \leq 1$ for all positive $\lambda$. It is equivalent to prove that for all positive $\lambda$
\begin{equation} \label{eqn:fun}
\frac{1}{t} \lambda^{\frac{1-t}{2t}} + \left( 1 - \frac{1}{t} \right) \lambda^{\frac{1}{2t}} \leq 1.
\end{equation}
Set $\displaystyle f(\lambda) := \frac{1}{t} \lambda^{\frac{1-t}{2t}} + \left( 1 - \frac{1}{t} \right) \lambda^{\frac{1}{2t}} - 1$. Let $a = \frac 1{2t}$. Then $a \in [1, \infty)$ and
\[ f(\lambda) = 2a \lambda^{a - \frac{1}{2}} + (1 - 2a) \lambda^a - 1.
\]
Since $f'(\lambda) = 2a (a - \frac{1}{2}) \lambda^{a - \frac{3}{2}} + a(1-2a) \lambda^{a-1} = a (2a - 1) \lambda^{a-1} (\lambda^{-1/2} - 1)$, $f(\lambda)$ attains the maximum value at $\lambda = 1$ for $\lambda > 0$. So \eqref{eqn:fun} holds since $f(\lambda) \leq f(1) = 0$ for $\lambda > 0$. Thus, by~\eqref{eqn:eigspe}, $\lambda_1(A \natural_t B) \leq 1$, that is, $A \natural_t B \leq I$.

For the case $t \in [\frac{1}{2}, 1)$, note that $1-t \in (0, \frac{1}{2}]$. By Lemma \ref{L:s-geo} (1) and Lemma \ref{L:Wass} (1),
if $ A \diamond_t B = B \diamond_{1-t} A \leq I$, then $A \natural_t B  = B \natural_{1-t} A \leq I$.
This completes the proof.
\end{proof}


\begin{corollary}
Let $A, B \in \mathbb{P}_m$ and $t \in [0,1]$. Then
\[ \| A \natural_t B \| \leq \| A \diamond_t B \|
\]
where $\| \cdot \|$ denotes the operator norm.
\end{corollary}

\begin{proof}
Let $c := \| A \diamond_t B \|$. Then $A \diamond_{t} B \leq cI$, and by Lemma \ref{L:Wass} (4)
\begin{displaymath}
\left( \frac{1}{c} A \right) \diamond_{t} \left( \frac{1}{c} B \right) = \frac{1}{c} (A \diamond_{t} B) \leq I.
\end{displaymath}
By Proposition \ref{P:sp-Wass} $\displaystyle \left( \frac{1}{c} A \right) \natural_{t} \left( \frac{1}{c} B \right) \leq I$. By Lemma \ref{L:s-geo} (4), $A \natural_{t} B \leq c I$, which yields the conclusion.
\end{proof}

To prove the weak log-majorization between two means, we use the property of the $k$th antisymmetric tensor power stated as below.

\begin{lemma}\cite[Problem I.6.12.]{Bh97}\label{lem:compound_monotone}
Let $A, B\in \mathbb{P}_m$ and $A\geq B$. Then
\[\Lambda^k A \geq \Lambda^k B.
\]
\end{lemma}

Now we are ready to prove the weak log-majorization between the spectral geometric mean and the Wasserstein mean.

\begin{theorem} \label{thm:main}
Let $A, B \in \mathbb{P}_m$ and $t \in [0,1]$. Then
\begin{equation} \label{eqn:speWas}
A \natural_t B \prec_{w \log} A \diamond_{t} B.
\end{equation}
\end{theorem}

\begin{proof}
Assume that $\Lambda^k (A\diamond_t B) \leq I$.
By \eqref{eqn:Wass-v}, we know that
\[A\diamond_t B = [I\nabla_t(A^{-1}\# B)]A[I\nabla_t(A^{-1}\# B)].
\]
Then taking the congruence transformation by $ \Lambda^k [I\nabla_t(A^{-1}\# B)]^{-1}$, we have
\[\begin{split}
\Lambda^k \{ [ I \nabla_t (A^{-1} \# B) ] A [ I \nabla_t (A^{-1} \# B) ] \} & \leq I \\
\Lambda^k [ I \nabla_t (A^{-1} \# B) ]^{-2} & \geq \Lambda^k A \\
\Lambda^k [ I \nabla_t (A^{-1}\# B) ]^{2} & \leq \Lambda^k A^{-1}.
\end{split}
\]
By Lemma~\ref{L:m-geo} (8) and Lemma~\ref{lem:compound_monotone},
\[ \begin{split}
\lambda_1(\Lambda^k [I \nabla_t (A^{-1} \# B)]^{2}) & = \lambda_1^2(\Lambda^k [I \nabla_t (A^{-1} \# B)]) \\
& \geq \lambda_1^2(\Lambda^k [I \#_t (A^{-1} \# B)]) = \lambda_1^2(\Lambda^k (A^{-1} \# B)^t ).
\end{split}
\]
So
\[ \lambda_1^{-1}(\Lambda^k A) \geq \lambda_1(\Lambda^k [I \nabla_t (A^{-1} \# B)]^{2}) \geq \lambda_1(\Lambda^k (A^{-1} \# B)^{2t} ).
\]
Thus, we have
\[\lambda_1(\Lambda^k (A \natural_t B)) = \lambda_1(\Lambda^k (A^{-1} \# B)^{2t} \Lambda^k A) \leq \lambda_1(\Lambda^k (A^{-1} \# B)^{2t}) \lambda_1 (\Lambda^k A) = 1,\]
where the first equality holds because $\lambda_{1}(A B) = \lambda_{1}(B A)$ and the inequality is valid since $(A^{-1}\# B)^{2t}, A \in \mathbb{P}_m$.
\end{proof}

\begin{remark}
By Theorem \ref{thm:main} together with known results appeared in Section \ref{Sec:properties} we have the following relation among matrix means:
\begin{displaymath}
A \#_{t} B \prec_{\log} \exp({(1-t) \log A + t \log B}) \prec_{\log} A \natural_{t} B \prec_{w \log} A \diamond_{t} B.
\end{displaymath}
Furthermore, \eqref{eqn:speWas} cannot be a log-majorization because of the determinantal inequality of the Wasserstein mean in \cite[Proposition 2.3]{HK19}:
\begin{displaymath}
\det \Omega(\omega; A_{1}, \dots, A_{n}) \geq \prod_{j=1}^{n} (\det A_{j})^{w_{j}},
\end{displaymath}
and equality holds if and only if $A_{1} = \cdots = A_{n}$, where $A_{j} \in \mathbb{P}_{m}$ for all $j = 1, \dots, n$. Also see Lemma \ref{L:s-geo} (6) and Lemma \ref{L:Wass} (6).
\end{remark}


\section{Inequalities of Wasserstein mean} \label{Sec:Was}

In this section, we explore some inequalities and the monotone convergence of Wasserstein mean with respect to the weak log-majorization. We finally propose a conjecture.

\begin{lemma} \label{L:supplement}
Let $A, B \in \mathbb{P}_{m}$. Then $A \diamond_{t} B \leq I$ for any $t \in (0,1)$ implies the following
\begin{itemize}
  \item[(i)] $A \leq B^{-1}$,
  \item[(ii)] $A \leq I$.
\end{itemize}
\end{lemma}

\begin{proof}
Assume that $A \diamond_{t} B \leq I$ for any $t \in (0,1)$.

For (i), by \eqref{eqn:speWas} $A \natural_{t} B \leq I$, which is equivalent to $A \leq B^{-1}$. See \cite[Remark 3.5]{HK22} and its references.

For (ii), by the definition of the Wasserstein mean~\eqref{eqn:Wass},
 \[A^{-1/2} [(1-t) A + t (A^{1/2} B A^{1/2})^{1/2}]^{2} A^{-1/2} \leq I.\]
  Taking congruence transformation by $A^{1/2}$ and using the operator monotonicity of square root map yield
\begin{displaymath}
(1-t) A + t (A^{1/2} B A^{1/2})^{1/2} \leq A^{1/2}.
\end{displaymath}
Taking congruence transformation by $A^{-1/2}$ implies
\begin{displaymath}
(1-t) I + t A^{-1} \# B \leq A^{-1/2}.
\end{displaymath}
By (i) $A^{-1} \geq B$, and thus, by the weighted arithmetic-geometric mean inequality
\begin{displaymath}
A^{-1/2} \geq (1-t) I + t A^{-1} \# B \geq (1-t) I + t B \geq I \#_{t} B = B^{t}.
\end{displaymath}
This implies that $A^{-1/2} \geq B^{t_{k}}$ holds for a sequence $t_{k} \in (0,1)$ converging to $0$, and hence, $A^{-1/2} \geq I$. That is, $A \leq I$.
\end{proof}


\begin{proposition} \label{P:weak-log-majorization}
Let $A, B \in \mathbb{P}_{m}$. If $A^{2} \diamond_{t} B^{2}\leq I$, then $\left( A \diamond_{t} B \right)^ 2 \leq I$
for any $\displaystyle t \in \left(0, \alpha_{B}/(\alpha_{B} + \beta_{B}) \right] \cup \left[ \beta_{A}/(\alpha_{A} + \beta_{A}), 1 \right)$, where $\alpha_{A} = \lambda_{m}(A), \beta_{A} = \lambda_{1}(A)$ and $\alpha_{B} = \lambda_{m}(B), \beta_{B} = \lambda_{1}(B)$.
\end{proposition}


\begin{proof}
It is enough to show that $A^{2} \diamond_{t} B^{2} \leq I$ implies $A \diamond_{t} B \leq I$.
Assume that $A^{2} \diamond_{t} B^{2} \leq I$.

\begin{itemize}
\item[Step 1.] For the case $\beta_{B} \leq 1$, equivalently $B \leq I$, we easily see that $A \diamond_{t} B \leq I$ from Lemma \ref{L:supplement} (ii) and \cite[Lemma 2.4]{KL20}.

\item[Step 2.] We prove that $A \diamond_{t} B \leq I$ when $B > I$. By Lemma \ref{L:supplement} (i) $A^{2} \leq B^{-2}$, and $A \leq B^{-1}$ by the operator monotonicity of square root map. So by \eqref{E:A-Wass}
\begin{displaymath}
A \diamond_{t} B \leq (1-t) A + t B \leq (1-t) B^{-1} + t B.
\end{displaymath}
All eigenvalues of $(1-t) B^{-1} + t B$ are of the form
\begin{displaymath}
f(t) := (1-t) \lambda^{-1} + t \lambda,
\end{displaymath}
where $\lambda \in (1, \beta_{B}]$ denotes an arbitrary eigenvalue of $B$. One can see that $f(t) \leq 1$ when $t \leq \frac{1}{1 + \lambda}$. Since $\lambda \leq \beta_{B}$, we have $\frac{1}{1 + \beta_{B}} \leq \frac{1}{1 + \lambda}$. Thus, $A \diamond_{t} B \leq I$ for all $\displaystyle t \in \left(0, \frac{1}{1 + \beta_{B}} \right]$.

\item[Step 3.] For the case $\alpha_{B} \leq 1 \leq \beta_{B}$, the preceding arguments in Step 2 yield that
\begin{displaymath}
\left( \frac{1}{\alpha_{B}} A \right) \diamond_{t} \left( \frac{1}{\alpha_{B}} B \right) \leq I
\end{displaymath}
for $\displaystyle t \in \left(0, \frac{1}{1 + \frac{\beta_{B}}{\alpha_{B}}} \right]$ because $\frac{1}{\alpha_{B}} B \geq I$ and the maximum eigenvalue of $\frac{1}{\alpha_{B}} B$ is $\frac{\beta_{B}}{\alpha_{B}}$. By the homogeneous property of two-variable Wasserstein mean in Lemma \ref{L:Wass} (4),
\begin{displaymath}
A \diamond_{t} B \leq \alpha_{B} I \leq I
\end{displaymath}
for all $\displaystyle t \in \left(0, \frac{\alpha_{B}}{\alpha_{B} + \beta_{B}} \right]$. So $(A \diamond_{t} B)^2 \leq I$ holds for $\displaystyle t \in \left(0, \frac{\alpha_{B}}{\alpha_{B} + \beta_{B}} \right]$.

\item[Step 4.] For $\displaystyle t \in \left[ \frac{\beta_{A}}{\alpha_{A} + \beta_{A}}, 1 \right)$, we have by Step 3 and the symmetric property of two-variable Wasserstein mean in Lemma \ref{L:Wass} (1)
\begin{displaymath}
\left( A \diamond_{t} B \right)^{2} = \left( B \diamond_{1-t} A \right)^{2} \prec_{w \log} B^{2} \diamond_{1-t} A^{2} = A^{2} \diamond_{t} B^{2}.
\end{displaymath}
because $\displaystyle 1-t \in \left(0, \frac{\alpha_{A}}{\alpha_{A} + \beta_{A}} \right]$.
\end{itemize}
These four steps complete the proof.
\end{proof}

Note that the interval
\[ \displaystyle \left(0, \frac{\alpha_{B}}{\alpha_{B} + \beta_{B}} \right] \cup \left[ \frac{\beta_{A}}{\alpha_{A} + \beta_{A}}, 1 \right)
\]
for $\alpha_{A} = \lambda_{m}(A), \beta_{A} = \lambda_{1}(A)$ and $\alpha_{B} = \lambda_{m}(B), \beta_{B} = \lambda_{1}(B)$, appeared in Proposition \ref{P:weak-log-majorization}, does not contain $1/2$. So it is an interesting question that it holds for all $t \in (0,1)$.

\begin{theorem} \label{T:square-Wass}
Let $A, B \in \mathbb{P}_{m}$ such that $(A^{-1} \# B)^{2} \leq A^{-2} \# B^{2}$ and $t \in [0,1]$. Then
\begin{displaymath}
\left( A \diamond_{t} B \right)^{2} \prec_{w \log} A^{2} \diamond_{t} B^{2}.
\end{displaymath}
\end{theorem}

\begin{proof}
Let $A, B \in \mathbb{P}_{m}$ such that $(A^{-1} \# B)^{2} \leq A^{-2} \# B^{2}$, and let $t \in [0,1]$.
Note that
\begin{displaymath}
A^{2} \diamond_{t} B^{2} = [I \nabla_{t} (A^{-2} \# B^{2})] A^{2} [I \nabla_{t} (A^{-2} \# B^{2})].
\end{displaymath}

Assume that $\Lambda^{k} (A^{2} \diamond_{t} B^{2}) \leq I$. Then $\Lambda^{k} A^{2} \leq \Lambda^{k} [I \nabla_{t} (A^{-2} \# B^{2})]^{-2}$, and by Theorem~\ref{thm:LH-ineq}
\begin{equation} \label{E:supplement}
\Lambda^{k} A \leq \Lambda^{k} [I \nabla_{t} (A^{-2} \# B^{2})]^{-1}.
\end{equation}
Since the square map is operator convex on $\mathbb{P}_{m}$, we have from the assumption
\begin{displaymath}
[I \nabla_{t} (A^{-1} \# B)]^{2} \leq I \nabla_{t} (A^{-1} \# B)^{2} \leq I \nabla_{t} (A^{-2} \# B^{2}).
\end{displaymath}
Since the $k$th compound is monotone on $\mathbb{P}_{m}$,
\begin{displaymath}
\Lambda^{k} [I \nabla_{t} (A^{-1} \# B)]^{2} \leq \Lambda^{k} [I \nabla_{t} (A^{-2} \# B^{2})],
\end{displaymath}
and hence, \eqref{E:supplement} yields
\begin{displaymath}
\Lambda^{k} A \leq \Lambda^{k} [I \nabla_{t} (A^{-2} \# B^{2})]^{-1} \leq \Lambda^{k} [I \nabla_{t} (A^{-1} \# B)]^{-2}.
\end{displaymath}
Thus, taking the congruence transformation by $\Lambda^{k} [I \nabla_{t} (A^{-1} \# B)]$ we obtain
\begin{displaymath}
\Lambda^{k} (A \diamond_{t} B) = \Lambda^{k} [I \nabla_{t} (A^{-1} \# B)] \left( \Lambda^{k} A \right) \Lambda^{k} [I \nabla_{t} (A^{-1} \# B)] \leq I,
\end{displaymath}
which completes the proof.
\end{proof}

\begin{remark}
The assumption in Theorem \ref{T:square-Wass} implies $(A^{-1/2} \# B^{1/2})^{2} \prec_{\log} A^{-1} \# B$. On the other hand, $A^{-1} \# B \prec_{\log} (A^{-1/2} \# B^{1/2})^{2}$ from \cite{AH94}. It implies that the spectra of $(A^{-1/2} \# B^{1/2})^{2}$ and $A^{-1} \# B$ are the same, which means that $(A^{-1/2} \# B^{1/2})^{2}$ and $A^{-1} \# B$ are unitarily similar.
\end{remark}

\begin{remark}
The special case of the assumption in Theorem \ref{T:square-Wass} is $(A^{-1} \# B)^{2} = A^{-2} \# B^{2}$, which is $A^{-1} \# B = (A^{-2} \# B^{2})^{1/2}$. By \cite[Theorem 3.1]{Hi94} it is equivalent to $A^{-1} B = B A^{-1}$, that is, $A B = B A$.
\end{remark}

Let $A, B \in \mathbb{P}_m$ and $t \in [0,1]$. Here, the notations $\nearrow_{\prec_{\log}}$ and $\searrow_{\prec_{\log}}$ mean to converge increasingly and decreasingly with respect to $\prec_{\log}$, respectively.
Hiai and Petz \cite{HP93} showed that the limit of the metric geometric mean for $A, B \in \mathbb{P}_m$ is the log-Euclidean mean, when $p$ goes to $0$:
\[ \lim_{p \to 0} (A^p \#_t B^p)^{1/p} = \exp((1-t) \log A + t \log B).
\]
Ando and Hiai \cite{AH94} gave the monotonicity of the metric geometric mean with respect to the log-majorization relation as follows:
\[ (A^p \#_t B^p)^{1/p} \prec_{\log} (A^q \#_t B^q)^{1/q}, \quad 0 < q \leq p.
\]
It is straightforward to verify
\[ (A^p \#_t B^p)^{1/p} \nearrow_{\prec_{\log}} \exp((1-t) \log A + t \log B) \quad \text{as} \quad p \searrow 0.
\]
Ahn, Kim and Lim \cite{AKL07} provided that the limit of the spectral geometric mean for $A, B \in \mathbb{P}_m$ is the log-Euclidean mean, when $p$ goes to $0$:
\[ \lim_{p \to 0} (A^p \natural_t B^p)^{1/p} = \exp((1-t) \log A + t \log B).
\]
Gan and Tam \cite{GT22} proved the monotonicity of the spectral geometric mean with respect to the log-majorization:
\[ (A^p \natural_t B^p)^{1/p} \searrow_{\prec_{\log}} \exp((1-t) \log A + t \log B) \quad \text{as} \quad p \searrow 0.
\]
The limit of the Wasserstein mean for $A, B \in \mathbb{P}_m$ is also the log-Euclidean mean, when $p$ goes to $0$ \cite{HK19}:
\[ \lim_{p \to 0} (A^p \diamond_t B^p)^{1/p} = \exp((1-t) \log A + t \log B).
\]
\begin{theorem} \label{thm:dec-Wass}
Let $A, B \in \mathbb{P}_{m}$ such that for any $0 < p \leq q \leq 1$
\begin{equation} \label{E:assumption}
(A^{-p} \# B^{p})^{1/p} \leq (A^{-q} \# B^{q})^{1/q}.
\end{equation}
Then
\begin{equation} \label{E:dec-Wass}
\left( A^{p} \diamond_{t} B^{p} \right)^{1/p} \searrow_{\prec_{w \log}} \exp ((1-t) \log A + t \log B) \quad \text{as} \quad p \searrow 0.
\end{equation}
\end{theorem}

\begin{proof}
The assumption with $p = 1/2$ and $q = 1$ yields $(A^{-1/2} \# B^{1/2})^{2} \leq A^{-1} \# B = (A^{-1/2})^{2} \# (B^{1/2})^{2}$.
Applying Theorem \ref{T:square-Wass} we have
\begin{displaymath}
\left( A^{1/2} \diamond_{t} B^{1/2} \right)^{2} \prec_{w \log} A \diamond_{t} B.
\end{displaymath}
The assumption with $p = 1/4$ and $q = 1/2$ yields $(A^{-1/4} \# B^{1/4})^{4} \leq (A^{-1/2} \# B^{1/2})^{2}$. By the Loewner-Heinz inequality in Theorem \ref{thm:LH-ineq} $(A^{-1/4} \# B^{1/4})^{2} \leq A^{-1/2} \# B^{1/2} = (A^{-1/4})^{2} \# (B^{1/4})^{2}$. Applying Theorem \ref{T:square-Wass} again we have
\begin{displaymath}
\left( A^{1/4} \diamond_{t} B^{1/4} \right)^{2} \prec_{w \log} A^{1/2} \diamond_{t} B^{1/2}.
\end{displaymath}
This implies
\begin{displaymath}
\left( A^{1/4} \diamond_{t} B^{1/4} \right)^{4} \prec_{w \log} (A^{1/2} \diamond_{t} B^{1/2})^{2} \prec_{w \log} A \diamond_{t} B.
\end{displaymath}
By induction we obtain
\begin{displaymath}
\left( A^{1/2^{k}} \diamond_{t} B^{1/2^{k}} \right)^{2^{k}} \prec_{w \log} A \diamond_{t} B.
\end{displaymath}
Since $\displaystyle \left( A^{1/2^{k}} \diamond_{t} B^{1/2^{k}} \right)^{2^{k}} \to \exp (1-t) \log A + t \log B)$ as $k \to \infty$, we conclude the desired property.
\end{proof}

One can ask whether \eqref{E:dec-Wass} holds without the assumption \eqref{E:assumption}, as the monotone convergence of metric and spectral geometric means with respect to the log-majorization. We close this paper by the following conjecture.

\begin{conjecture}
Let $A, B \in \mathbb{P}_{m}$ and $t \in (0,1)$. Then
\begin{displaymath}
\left( A^{p} \diamond_{t} B^{p} \right)^{1/p} \searrow_{\prec_{w \log}} \exp ((1-t) \log A + t \log B) \quad \text{as} \quad p \searrow 0.
\end{displaymath}
\end{conjecture}

\noindent  \textbf{Acknowledgement}

\noindent
We deeply thank to anonymous reviewer for valuable comments. The work of S. Kim was supported by the National Research Foundation of Korea grant funded by the Korea government (MSIT) (No. NRF-2022R1A2C4001306).

\bibliographystyle{abbrv}

\bibliography{refs}

\end{document}